\documentstyle[12pt,amscd,amssymb]{article}

\newtheorem{theorem}{Theorem}[section]
\newtheorem{lemma}[theorem]{Lemma}
\newtheorem{proposition}[theorem]{Proposition}

\newtheorem{corollary}[theorem]{Corollary}

\newcommand{\pa}{P_{A \times A'}}
\newcommand{\pbb}{P_{B \times B'}}
\newcommand{\cA}{{\cal A}}

\newcommand{\tim}{^\times}
\newcommand{\an}{^{an}}
\newcommand{\dan}{^{' an}}

\newcommand{\inv}{^{-1}}
\newcommand{\iso}{\stackrel{\sim}{\longrightarrow}}

\newcommand{\hdra}{H^1_{dR}(A)}
\newcommand{\hdrb}{H^1_{dR}(B)}
\newcommand{\hdrg}{H^1_{dR}(G)}
\newcommand{\hdrt}{H^1_{dR}(T)}

\newcommand{\ha}{H^1(A,{\cal{O}})}
\newcommand{\hb}{H^1(B,{\cal{O}})}

\newcommand{\lie}{\mbox{\rm Lie}}
\newcommand{\indif}{\mbox{\rm Inv}}
\newcommand{\ext}{\mbox{\rm Ext}}
\newcommand{\pic}{\mbox{\rm Pic}}

\newcommand{\Inf}{\mbox{\rm Inf}}
\newcommand{\lra}{\longrightarrow}
\newcommand{\lla}{\longleftarrow}

\newcommand{\mathq}{{\mathbb{Q}}}
\newcommand{\mathg}{{\mathbb{G}}}

\newcommand{\pair}{(\; \; , \; \;)}

\newcommand{\cL}{{\cal L}}

\newcommand{\cO}{{\cal O}}
\newcommand{\eps}{{K[\epsilon]}}
\newcommand{\ov}{\overline}
\newcommand{\cX}{{\frak X}}
\newcommand{\hX}{\hat{\frak X}}
\newcommand{\barX}{\overline{\frak X}}
\newcommand{\PP}{({\Bbb P}_K^1)^t}
\newcommand{\bG}{{\Bbb G}}
\newcommand{\bP}{{\Bbb P}}
\newcommand{\Q}{{\Bbb Q}}

\title{$p$-adic height pairings on abelian varieties with semistable ordinary 
reduction}
\author{Adrian Iovita \and  Annette Werner }
\date{}

\begin{document}
\maketitle
\centerline{\bf Abstract}

We prove that for abelian varieties with semistable ordinary reduction the 
$p$-adic Mazur-Tate height pairing is induced by the unit root splitting of the 
Hodge filtration on the first deRham cohomology.

\section*{Introduction}
The main goal of this paper is to prove that the $p$-adic Mazur-Tate height 
pairing on an abelian variety with semistable ordinary reduction is induced by the 
unit root splitting of the Hodge filtration on its first deRham cohomology. 

$p$-adic height pairings are the $\Q_p$-valued counterparts of the real-valued 
N\'eron-Tate height pairings on abelian varieties over a global field $F$. As the 
N\'eron-Tate pairings they can be decomposed into local contributions, one for 
each finite place of the ground field $F$. At the places not dividing $p$, these 
local contributions are  basically given by the local N\'eron heights. 
Hence only at the places over $p$ something genuinely $p$-adic  happens.

If the abelian variety $A$ has semistable ordinary reduction at a place over $p$, 
there are (at least) three interesting candidates for such a local  $p$-adic height pairing on 
$A$. First of all, there is the Mazur-Tate height defined with splittings of the 
Poincar\'e biextension. Then there is Schneider's ``norm-adapted'' $p$-adic 
height, and finally one can use the unit root splitting of the Hodge filtration
 on
the first deRham cohomology to define a $p$-adic height pairing. 
In the first section of this paper we explain some details of these constructions.
In particular we show how to pass from a pairing defined  by a splitting the 
Hodge filtration of the first deRham cohomology to a pairing in the Mazur-Tate
style.

It is well-known that if $A$ has {\em good} ordinary reduction, the Mazur-Tate 
height pairing coincides with Schneider's height pairing, cf. \cite{mata}. Besides, 
Coleman has proven in \cite{co} that in this case  the Mazur-Tate pairing  is 
given by the unit root splitting. Hence in the good ordinary reduction case all 
three definitions give the same pairing.

What happens in the semistable ordinary reduction case? In \cite{we} it is 
shown that the Mazur-Tate pairing differs in general from Schneider's height 
pairing. Hence the question remains if the height pairing given by the unit root 
splitting coincides with one of those pairings. We answer this question by showing 
in Theorem 3.6 that it is equal to the Mazur-Tate pairing. 

In the course of the proof we use the Raynaud extension
 which is a lift of the semiabelian variety in the 
reduction of $A$ to a semiabelian rigid analytic variety lying over $A$. Its 
abelian quotient $B$ has good reduction. By Coleman's result,  the Mazur-Tate 
height on $B$ is given by the unit root splitting. A result of the second author
 shows how the  
Mazur-Tate height pairings on $A$ and $B$ are related. Using results of  LeStum,
 the deRham cohomologies of $A$ and $B$ 
are connected by a diagram involving the deRham cohomology of the Raynaud 
extension. Besides, Coleman and the first author have shown how the Frobenius on the deRham 
cohomologies can be described explicitely. We use these facts in section 2
to relate the unit 
root splittings for $A$ and $B$.  

However, the previous results cannot simply be put together to prove the desired
Theorem 3.6. The reason is that the step from height pairings defined with splittings
of the Hodge filtration to Mazur-Tate height pairings involves the universal
vectorial extension of the abelian variety - and we have no result relating
the universal vectorial extensions of $A$ and $B$. 
Hence we go the other way and start with the Mazur-Tate pairing on $A$.
We prove in section 3 that it is of an analytic nature. To be precise,
we show that it gives rise to a $p$-adic analytic splitting of the rational 
points
of the universal vectorial 
extension of $A$. Hence the Mazur-Tate height is in fact 
induced by a splitting of the
Hodge filtration on the first deRham cohomology, and we can adapt some ideas from
\cite{mame} to deduce that this is the unit root splitting. 

In an appendix we study  the Hodge filtration of the first deRham cohomology 
of a semiabelian variety and show that it is given by its invariant 
differentials. This answers a question raised in \cite{ls2}.

{\bf Acknowledgements: }We are very grateful to Christopher Deninger for his
interest in these results. Much of the research  on this paper was
carried out while the first author was a guest of the SFB 478 in M\"unster. 
He is grateful to this institution for its
support and hospitality during two visits. The first author was
also partially supported by an NSF research grant. 

\section{$p$-adic height pairings}

In this section we collect various facts about $p$-adic height pairings, thereby 
fixing our notation.

Let $A_F$ be an abelian variety over a number field $F$, and let $A'_F$ denote its 
dual abelian variety. 
 
Recall that the classical N\'eron-Tate height pairing on $A_F$ 
\[\pair_N: A'_F(F) \times A_F(F) \rightarrow {\mathbb{R}}\]
can be decomposed in a sum of local pairings
\[\pair_{N,v}: ({\mbox{Div}}^0(A_{F_v}) \times Z^0(A_{F_v}/F_v))' \rightarrow 
{\mathbb{R}},\]
where $v$ runs over the places of $F$, and $F_v$ denotes the $v$-adic completion 
of $F$. 
Here $({\mbox{Div}}^0(A_{F_v}) \times Z^0(A_{F_v}/F_v))'$ is the set of all pairs
$(D,z)$ such that $D$ is a divisor on $A_{F_v}$, algebraically equivalent to zero, 
and $z$ is a cycle $\sum_i n_i a_i$ with $a_i \in A_{F_v}(F_v)$ and $\sum_i n_i 
=0$, such that the supports of $D$ and $z$ are disjoint.
These local pairings can be characterized by a list of axiomatic properties, see 
\cite{ne}. For all rational functions $f$ on $A_{F_v}$ we have 
$(\mbox{div}(f), \sum n_i a_i)_{N,v} = \sum n_i \log|f(a_i)|_v$. 

Now we want to investigate $p$-adic height pairings, i.e. pairings 
\[\pair_p: A'_F(F) \times A_F(F) \rightarrow \mathq_p.\]
The initial datum is a collection of continuous homomorphisms $\rho_v: F_v\tim 
\rightarrow \mathq_p$, one for each finite place $v$ of $F$, such that $\sum_v  
\rho_v = 0$. (This substitutes the collection $(\log|\,|_v)$ from the classical 
case.)
For example, we could take 
\[ \rho_v = \left\{ \begin{array}{ll}
   \log_p \circ N_{F_v/ {\mathbb{Q}}_p} & \mbox{ if $v|\,p$}\\
   \log_p \circ |\;\; |_l \circ N_{F_v/{\mathbb{Q}}_l} & \mbox{ if $v | \! \!  / 
\,p$,}
\end{array} \right. \]
where $\log_p: {\mathbb{Q}}^\times_p \rightarrow \Q_p$ is the branch of the 
$p$-adic logarithm vanishing on $p$.

Note that for all places $v$ not dividing $p$, any continuous homomorphism 
$\rho_v$ is unramified, i.e. it maps the elements of absolute value one to zero. 
Hence it is proportional to the valuation map $v$, and also to $\log|\,|_v$.

Assume for a moment that we can construct a collection of bilinear maps 
\[\pair_v: ({\mbox{Div}}^0(A_{F_v}) \times Z^0(A_{F_v}/F_v))' \rightarrow 
\mathq_p,\]
one for each finite place of $F$,
with the following properties:

1) $({\rm div}f, \sum n_i a_i )_v = \sum n_i \rho_v(f(a_i))$ for all rational 
functions $f$ on $A_{F_v}$. 

2) $(t_a^\ast D, t_a^\ast z)_v = (D,z)_v$ for all $a\in A_{F_v}(F_v)$ and all 
pairs $(D,z)$.

3) If $\rho_v$ is unramified, i.e. $\rho_v(x) = c v(x)$ for some constant $c$, 
then 
$\pair_v$ is proportional to the N\'eron pairing.

With these data we can construct a global $p$-adic height pairing as follows: For 
$a \in A_F(F)$ and $a' \in A'_F(F)$ choose a divisor $D$, algebraically equivalent 
to zero, whose class is $a'$, and whose support contains neither $a$ nor the unit 
element $ 0 \in A_F(F)$. Then put
\[(a',a) = \sum_{v}(D, a-0)_v \in \mathq_p,\]
where we lift $D$ and the cycle $a-0$ to the variety $A_{F_v}$. This is 
well-defined, i.e. the sum is finite, and independent of the choice of $D$. Using 
property 2) of the local height pairings one can show that it is bilinear.

One recipe to construct local pairings with the desired properties was given by 
Mazur and Tate in \cite{mata}. 

Let us fix a place $v$ and a continuous homomorphism $\rho_v: F_v\tim \rightarrow 
\mathq_p$. From now on we will work in a local setting. To save subscripts, let us 
put $\rho= \rho_v$ and $K = F_v$. Besides, denote by $R$ the ring of integers in 
$K$ and by $k$ the residue field. Recall that we call $\rho$ unramified, if it 
maps $R\tim$ to zero. 

Let $\cA$ respectively $\cA'$ be N\'eron models of $A = A_F \times_F K$ 
respectively $A' = A'_{F} \times_F K$ over $R$. Besides, let $P= P_{A \times A'}$ 
be the Poincar\'e biextension expressing the duality between $A$ and $A'$, see 
\cite{sga7},  VIII, 1.4. 
$P$ induces an isomorphism $A' \simeq \underline{\ext}^1(A, \ \mathg_{m})$ of 
sheaves in the fppf-site over $K$, mapping a point $a' \in A'(T)$ to the 
restriction $P_{A_T \times \{a'\}}$. We will often identify points in $A'$ with 
extensions using this isomorphism.

We call a map
$\sigma: P(K) \rightarrow \mathq_p$ a $\rho$-splitting if $\sigma$ behaves 
homomorphically with respect to both group laws and if $\sigma(\alpha x) = 
\rho(\alpha) + \sigma(x)$ for all $\alpha \in K\tim$ and $x \in P(K)$, see 
\cite{mata}, 1.4.

For any $\rho$-splitting $\sigma$ we can define a bilinear map $\pair_{\sigma}: 
({\mbox{Div}}^0(A) \times Z^0(A/K))' \rightarrow \mathq_p$ with the properties 1) 
and 2) as follows: For any divisor $D \in {\mbox{Div}}^0(A)$ let $d$ be the point 
in $A'(K)$ corresponding to its class. Then $P_{A \times \{d\}}$ is the extension 
of $A$ by ${\mathbb{G}}_{m}$ corresponding to the class of $D$, hence $D$ gives 
rise to a rational section $s_D$ of $P_{A \times \{d\}} \rightarrow A$, which is a 
morphism on the complement of the support of $D$.
Hence we can put
\[(D,\sum n_i a_i)_{\sigma} = \sum n_i \sigma(s_D(a_i))\]
to get our local pairing. Conversely, by \cite{mata}, 2.2, every local height 
pairing with 1) and 2) comes in fact from a uniquely defined $\rho$-splitting 
$\sigma$.

In \cite{mata} canonical  $\rho$-splittings are defined in several cases. We will 
discuss two of them. The first one is the case that $\rho$ is unramified. 

Note that $P$ can be extended to a biextension $P_{\cA^0 \times \cA'}$ of $\cA^0$ 
and $\cA'$ by $\mathg_{m,R}$, where $\cA^0$ is the identity component of $\cA$, 
see \cite{sga7}, VIII, 7.1.

Since $\rho$ is unramified, there exists a unique $\rho$-splitting $\sigma$ 
vanishing on $P_{\cA^0 \times \cA'}(R)$, see \cite{mata}, 1.9. If $\rho$ is 
equal to $\log|\,|_v$, then the pairing $\pair_{\sigma}$ coincides with the local 
N\'eron pairing $\pair_{N,v}$ at our fixed place $v$ by \cite{mata}, 2.3.1. Hence 
we find that for unramified $\rho$ the Mazur-Tate pairing $\pair_{\sigma}$ 
fulfills property 3) from above.

Note that this means that it remains to construct local pairings for ramified 
continuous homomorphisms $\rho: K\tim \rightarrow \mathq_p$, which will only 
involve places $v$ over $p$.

Another case where Mazur and Tate construct a canonical $\rho$-splitting is the 
case that $A$ has semistable ordinary reduction, i.e. that the formal completion 
of the special fibre $\cA_k$ along the zero section  is isomorphic to a product of 
copies of ${\mathbb{G}}_m^f$ over the algebraic closure of $k$. This is equivalent 
to the fact that $\cA_k^0$ is an extension of an ordinary abelian variety by a 
torus, see \cite{mata}, 1.1. Note that in particular $\cA$ has semistable 
reduction.
Denote by $\cA^f$ and $\cA'^f$ the formal completions along the zero sections in 
the special fibres, and by $P_{\cA^0 \times \cA'}^f$ the formal completion of 
$P_{\cA^0 \times \cA'}$ along the preimage of the zero section of $\cA^0_k \times 
\cA'_k$ under the projection map. Then $P_{\cA^0 \times \cA'}^f$ is a formal 
biextension of $\cA^f$ and $\cA'^f$ by $\mathg_m^\wedge$, the formal completion of 
$\mathg_m$ along the special fibre. 
If $A$ has semistable ordinary reduction, $P_{\cA^0 \times \cA'}^f$ admits a 
unique trivialization, hence there is a uniquely defined biextension splitting 
$\tilde{\sigma}: P_{\cA^0 \times \cA'}^f \rightarrow \mathg_m^\wedge$, see 
\cite{mata}, 5.11.2. Then there exists a unique $\rho$-splitting $\sigma: P(K) 
\rightarrow \mathq_p$ such that for all $x \in P_{\cA^0 \times \cA'}^f(R)$ we have 
$\sigma(x) = \rho(\tilde{\sigma}(x))$, see \cite{mata}, 1.9. We call the 
corresponding local height pairing the canonical Mazur-Tate pairing in the 
semistable ordinary reduction case. 

Note that if additionally $\rho$ is unramified, we get the same $\rho$-splitting 
as defined previously (see \cite{mata}, p.204). 

There is also another way of constructing in certain cases for ramified $\rho$ a 
$\rho$-splitting, namely the norm-adapted Schneider height, cf \cite{sch}. 
By \cite{mata}, 1.11.6 it coincides with the Mazur-Tate height if $A$ has good 
ordinary reduction. A formula for the difference of  
these two $p$-adic height pairings in the case of semistable ordinary 
reduction can be found in \cite{we}, 7.2.

Let us now give a brief account of another approach to construct local $p$-adic 
height pairings, which uses splittings of the Hodge filtration of the first deRham 
cohomology. (See e.g. \cite{cogr} or \cite{za}.)

Let $V'$ be the vector group corresponding to $e^* \Omega_{A/K}^1$, where $e$ is 
the unit section of $A$. By $I'$ we denote the universal vectorial extension of 
the dual abelian variety $A'$, cf. \cite{mame}, chapter 1, \S 1. Then $I'$ sits in an exact sequence
\[ 0 \longrightarrow V' \longrightarrow I' \longrightarrow A' \longrightarrow 0,\]
and it is universally repellent with respect to extensions of $A'$ by vector 
groups, see \cite{mame}. $I'$ represents the sheaf 
$\underline{\rm Extrig}_K(A, \mathg_{m})$ of rigidified extensions of $A$ by 
$\mathg_{m}$. This is the Zariski sheaf associated to the presheaf mapping $S$ to 
${\rm Extrig}_S(A_S, \mathg_{m,S})$, the group of isomorphism classes of pairs 
$(E,t)$, where $E$ is an extension of $A_S$ by $\mathg_{m,S}$ and $t$ is a section 
on the level of the first infinitesimal extensions. To be precise, $t$ is a 
morphism of $S$-pointed $S$-schemes $t : \Inf^1(A_S/S) \rightarrow E$ such that 
the diagram
\begin{eqnarray*}
\begin{CD}
E @>>> A_S \\
@A{t}AA @AAA \\
\Inf^1(A_S/S) @= \Inf^1(A_S/S)
\end{CD}
\end{eqnarray*}
commutes. Here $\Inf^1$ is the first infinitesimal 
neighbourhood of the unit section in the sense of \cite{ega4}, 16.1.2.  If we 
identify $A'$ with the sheaf $\underline{\rm Ext}^1_K(A, \mathg_{m})$, the 
projection $I' \rightarrow A'$ corresponds to the map "forget the rigidification". 
 
Moreover, by \cite{mame}, chapter 1, \S 4, there is an isomorphism $\mbox{Lie}\, I' 
\iso \hdra$ such that the following diagram commutes
\begin{eqnarray*}
\begin{CD}
0 @>>> \mbox{Lie}\,V' @>>> \mbox{Lie}\,I' @>>> \mbox{Lie}\,A' @>>> 0 \\
@. @VV{\simeq}V @VV{\simeq}V  @VV{\simeq}V @. \\
0 @>>> H^0(A, \Omega^1) @>>> \hdra @>>> H^1(A, {\cal O}) @>>> 0,
\end{CD}
\end{eqnarray*}
where the left vertical map is obvious and the right one is given by \cite{mu}, 
p.130, and where the lower horizontal sequence is given by the Hodge filtration on 
$\hdra$.

For any smooth commutative $K$-group $H$ the $K$-rational points $H(K)$ 
define a Lie group over $K$ (in the sense of \cite{bo}, chapter III) whose Lie 
algebra coincides with the algebraic Lie algebra ${\mbox{Lie}} \, H$.

Suppose now that $r: H^1(A, {\cal O}) \rightarrow \hdra$ is a splitting of the 
Hodge filtration. 
This induces a splitting $\mbox{Lie} \, A' \rightarrow \mbox{Lie} \, I'$ which can 
be lifted to a Lie group homomorphism in a suitable neighbourhood of the unit 
section by    \cite{bo}, III, \S 7, Theorem 3. It is easy to see that this can 
be extended to a (uniquely determined) splitting  $\eta: A'(K) \rightarrow I'(K)$ 
of the projection $I'(K) \rightarrow A'(K)$ with $\mbox{Lie}\, \eta =r$, cf. 
\cite{za}, Theorem 3.1.3. With other words, this is a splitting of the forgetful 
homomorphism 
\[\mbox{Extrig}(A, \mathg_{m}) \rightarrow \mbox{Ext}^1(A, \mathg_{m}),\]
i.e. we found a multiplicative way of associating to an extension $X$ of $A$ by 
$\mathg_{m}$ a rigidification. Note that a rigidification on $X$ is the same as a 
splitting of the corresponding sequence of Lie algebras
\[ 0 \longrightarrow \lie\, \mathg_{m} \longrightarrow \lie\, X \longrightarrow 
\lie\,  A \longrightarrow 0.\]

Now take a divisor $D \in {\mbox{Div}}^0(A)$ whose class gives a point 
$d \in A'_K(K)$. Then $P_{A \times \{d\}}$, the extension corresponding to $d$, is 
endowed with a rigidification, which in turn induces a splitting 
\[t_d: \mbox{Lie}\, P_{A \times \{d\}}(K) = \mbox{Lie}\, P_{A \times \{d\}} 
\rightarrow \mbox{Lie}\, \mathg_{m} = \mbox{Lie} \, K\tim \]
of the Lie algebra sequence corresponding to the extension $P_{A \times \{d\}}$.

We fix again a continuous, ramified  homomorphism $\rho: K\tim \rightarrow \Q_p$. 
By \cite{za}, p. 319  
we have $\rho = \delta \circ \lambda$, where
$\delta: K \rightarrow \mathq_p$ is a $\mathq_p$-linear map, and 
$\lambda: K\tim \rightarrow K$ is a branch of the $p$-adic logarithm. 
Using once more \cite{bo}, III, \S 7, Theorem 3 we find a uniquely determined 
homomorphism of Lie groups 
\[\gamma_d: P_{A \times \{d\}}(K)  \rightarrow K\]
such that $K\tim \rightarrow P_{A \times \{d\}}(K) 
\stackrel{\gamma_d}{\rightarrow}K$
is the homomorphism $\lambda$ and such that $\mbox{Lie}\, \gamma_d = \mbox{Lie}\, 
\lambda \circ t_d$, cf. \cite{za}, Theorem 3.1.7. These maps $\gamma_d$ fit together 
to a $\lambda$-splitting $\gamma$ of $P(K)$. Hence $\delta \circ \gamma$ is a 
$\rho$-splitting of $P(K)$. 

Now we can define a pairing
\[\pair_r: ({\mbox{Div}}^0(A) \times Z^0(A/K))' \rightarrow \mathq_p\]
by $(D,\sum n_i a_i)_r = \sum n_i \, \delta\circ \gamma(s_D(a_i))$, where as above $s_D$ is a 
rational section of $P_{A_{K} \times \{d\}}$ corresponding to $D$. This map is 
bilinear and has properties 1) and 2). Obviously, $\pair_r$ is just the 
height pairing associated to the $\rho$-splitting $\delta \circ\gamma$.

Hence every splitting of the Hodge filtration induces a $\rho$-splitting 
on $P(K)$ such that
the corresponding height pairings are the same. What about the other direction?

Suppose that  $\tau: P(K) \rightarrow K$ is a $\lambda$-splitting.
Then we can define a 
splitting
\[\eta(\tau): A'(K) \longrightarrow I'(K)\]
of the projection map by associating to every $a' \in A'(K)$ the extension
$P_{A \times \{a'\}}$ endowed with the rigidification $
(\lie \lambda)\inv \circ \lie (\tau_{|P_{A \times \{a'\}}(K)}): 
 \lie P_{A \times \{a'\}}(K) {\longrightarrow}
 \lie K \iso \lie K\tim$. 
If $\eta(\tau)$ is analytic (in the sense of \cite{bo2}), then it induces
a Lie algebra splitting $r : \lie A' \rightarrow \lie I'$, hence a splitting of the
Hodge filtration of $\hdra$. This construction is converse
to the previous association $r \mapsto \gamma$. We will see in Proposition 3.2 that
$\eta(\tau)$ is analytic if $\tau$ is an analytic map.

\section{Splittings of the Hodge filtration of $\hdra$}

Let, as in section 1, $K$ be a non-archimedean local field of characteristic $0$
 with ring of integers 
$R$ and residue field $k$. Besides, let $A$ be an abelian variety over $K$ with 
ordinary semistable reduction,  $A'$ 
its  dual abelian variety, and $\cA$ respectively $\cA'$ their N\'eron models. 
We assume that the torus parts
in the reductions of $\cA^0$ and $\cA'_0$ are split. For our purposes, this is no 
restriction, since 
height pairings are compatible with base change.
  
Let us recall some facts about the rigid analytic uniformization of $A$ and $A'$.
There is an extension of algebraic groups over $K$
\[0 \longrightarrow T\stackrel{g}{\longrightarrow} G \stackrel{p}{\longrightarrow} 
B \longrightarrow 0,\]
such that $T$ is a split torus of dimension $t$ over $K$, and $B$ is
 an abelian variety over $K$ with good reduction.
There is also  a rigid analytic homomorphism $\pi: G\an \rightarrow A\an$ inducing 
a short exact sequence of rigid analytic groups over $K$
\[0 \longrightarrow \Gamma^{an} \stackrel{i}{\longrightarrow} G^{an} 
\stackrel{\pi}{\longrightarrow} A\an \longrightarrow 0,\]
where $\Gamma$ is the constant group scheme corresponding to a free 
${\mathbb{Z}}$-module $\Gamma$ of rank $t$. (See \cite{bolu}, section 1, and 
\cite{ray}.) 

Let $\Gamma'$ be the character group of $T$. We denote the corresponding constant 
$K$-group scheme also by $\Gamma'$. Fix a dual abelian variety $(B', \pbb)$ of 
$B$, where $\pbb$ is the Poincar\'e biextension expressing the duality. Then $G$ 
corresponds to a homomorphism $\phi': \Gamma' \rightarrow B'$ (see e.g. 
\cite{sga7}, VIII, 3.7).

The embedding $i: \Gamma \rightarrow G$ induces a homomorphism $\phi: \Gamma 
\stackrel{i}{\rightarrow} G \stackrel{p}{\rightarrow}B$, which gives us an 
extension $G'$ (again by \cite{sga7}, VIII, 3.7)
\[0 \longrightarrow T' \longrightarrow G' \stackrel{p'}{\longrightarrow} B' 
\longrightarrow 0,\]
where $T'$ is the split torus of dimension $t$ over $K$ with character group 
$\Gamma$. There is a short exact sequence 
\[0 \longrightarrow \Gamma^{' an} \stackrel{i'}{\longrightarrow} G^{' an} 
\stackrel{\pi'}{\longrightarrow} A^{' an} \longrightarrow 0.\]

In \cite{ls1} and \cite{ls2} it is shown how the rigid analytic 
uniformization of $A$ can be used to describe the deRham cohomology of $A$.  
Let us first fix some notation. 

Let $X$ be a commutative $K$-group variety, hence smooth over $K$. 
The space of invariant differentials 
$\indif(X)$ of $X$ is the space of sections of $e^\ast \Omega^1_{X/K}$, where $e: 
\mbox{\rm Spec}(K) \rightarrow X$ is the unit section. 
Note that $\indif(X)$ can be identified with the space of global differentials
$\omega \in \Gamma(X, \Omega_X^1)$ satisfying $m^\ast \omega = p_1^\ast \omega
+ p_2^\ast \omega$, where $m,p_1,p_2: X \times X \rightarrow X$ denote
multiplication and projections, respectively.
Besides, all invariant differentials
are closed.

 If $Z$ is a rigid analytic $K$-group variety, we define $\indif(Z)$ in the same 
way, using the rigid analytic differentials $\Omega^1_{Z/K}$ as defined e.g. in 
\cite{bkkn}. There is a natural GAGA-isomorphism $\indif(X) \simeq 
\indif(X^{an})$. 

Let us  denote by $H^1_{dR}(X)$, respectively, $H^1_{dR}(Z)$ the first deRahm 
cohomology group of the algebraic variety $X$, respectively, the rigid analytic 
variety $Z$. Hence $H^1_{dR}(X)$ ist the first hypercohomology of the complex 
$(0 \rightarrow {\cal O}_{X} \rightarrow \Omega^1_{X/K} \rightarrow \Omega^2_{X/K} \rightarrow 
\ldots)$ and, similarly,  $H^1_{dR}(Z)$ is the first hypercohomology of the 
complex $(0 \rightarrow {\cal O}_{Z} \rightarrow \Omega^1_{Z/K} \rightarrow \Omega^2_{Z/K} 
\rightarrow \ldots )$. By \cite{ki}, there is a GAGA-isomorphism
\[ H^1_{dR}(X) \simeq H^1_{dR}(X^{an}),\]
which we tacitely use to identify these groups.

By \cite{ls1} we have a commutative diagram with exact rows and columns:
\begin{eqnarray*}
\begin{CD}
@. 0 @. 0 @. @.\\
@. @VVV @VVV @.@.\\
~~~ \quad 0 @>>> \indif(B) @>>> \hdrb @>{\delta}>>  H^1(B, {\cal{O}}_{B}) @>>> 0\\
@.  @VVV @VV{p^\ast}V   @V{\simeq}V{\beta}V @.\\
(1) \quad 0 @>>> \indif(G) @>>> \hdrg @>>> H @>>> 0 \\
@. @VVV @V{g^\ast}VV @.@.\\
@. \indif(T) @>>{\simeq}> \hdrt @.@.\\
@. @VVV @VVV @. @.\\
@.0 @. 0 @. @. 
\end{CD}
\end{eqnarray*}
where $H$ is by definition the vector space making the middle horizontal sequence 
exact. The first row of this diagram is again induced by the Hodge filtration on 
$\hdrb$. We show in the appendix that in fact the middle row is also induced by 
the Hodge filtration on $\hdrg$.

Besides, we have a commutative diagram with exact rows and columns 
for $\hdra$ (see for example \cite{coleman_iovita}.)

\begin{eqnarray*}
\begin{CD}
@. @. 0 @. 0 @. \\
@. @. @VVV @VVV @. \\
@.@. Hom(\Gamma,K) @>{=}>> Hom(\Gamma,K) @.\\
@. @. @VVV   @VVV @.\\
(2) \quad 0 @>>> \indif({A}) @>>> \hdra @>>> H^1(A,{\cal{O}}) @>>> 0\\
@. @V{\alpha}V{\simeq}V @VV{\pi^\ast}V @VV{\gamma}V @.\\
~~~ \quad 0 @>>> \indif({G}) @>>> \hdrg @>>> H @>>> 0\\
@. @.  @VVV  @VVV @. \\
@. @. 0 @. 0 @.
\end{CD}
\end{eqnarray*}
Here the middle row is again induced by the Hodge filtration on $\hdra$.

With these two diagrams we can lift any splitting $r: \hb \rightarrow \hdrb$ of 
the Hodge filtration of $\hdrb$ to a splitting $L(r) : \ha \rightarrow \hdra$ of 
the Hodge filtration of $\hdra$.
First lift $r$ to a splitting 
\[p^\ast \circ r \circ \beta\inv: H \rightarrow \hdrg\]
of the middle row in diagram $(1)$. Take the corresponding splitting
$s: \hdrg \rightarrow \indif(G)$ in the other direction, and lift this to
\[\alpha\inv \circ s \circ \pi^\ast: \hdra \rightarrow \indif(A)\]
with the help of diagram $(2)$. 

Then let $L(r): \ha \rightarrow \hdra$ be the splitting in the other direction. A 
diagram chase shows that $L(r)$ is the unique splitting of the Hodge filtration of 
$\hdra$ making the following diagram commutative:
\begin{eqnarray*}
\begin{CD}
 @. \hdra  @<{L(r)}<< \ha\\
 @. @VV{\pi^\ast}V @VV{\gamma}V \\
 @. \hdrg @. H\\
 @. @AA{p^\ast}A @VV{\beta\inv}V \\
@.\hdrb @<{r}<< {\hb}\\
\end{CD}
\end{eqnarray*}

For all $X\in \{A,G,T,B,\Gamma \}$, the $K$-vector space $H^1_{dR}(X)$
(which is by definition $Hom(\Gamma, K)$ if $X=\Gamma$) can be endowed 
with a ($K$-linear) Frobenius operator
$$
\varphi _X\colon H^1_{dR}(X)\rightarrow H^1_{dR}(X)
$$
(which is canonical for all $X$'s except for $X=A$, when it depends on the 
choice of a branch of log on $K^*$.)
By definition $\varphi _{\Gamma }=id$, $\varphi _T$= 
multiplication by $q=p^{[k:{\bf F}_p]}$. See \cite{coleman_iovita}
for properties of these maps.

Let now $X\in\{A,B,G\}$, then we put $H(X)=H$ if $X=G$ and $H(X)=H^1(X,\cO)$ if 
$X=B$ or $A$. For each $X\in\{A,G,B\}$ we define  the unit root subspace, $W_X
\subset H^1_{dR}(X)$ to be the subspace on which $\varphi _X$ acts with 
slope $0$.  
Let us  remark that although $\varphi_A$ 
depends on  the choice of a branch of the $p$-adic logarithm on $K^*$, 
$W_A$ is canonical (see \cite{illusie} or \cite{iovita}.)

So let us define $r_X$ to be the 
unique splitting $H(X) \rightarrow H^1_{dR}(X)$ such that  Im$(r_X)\subset
W_X$. By $s_X: H^1_{dR}(X) \rightarrow \indif(X)$ we denote the 
corresponding splitting in the other direction. Either of $s_X$ or $r_X$ will be 
called ``the unit root splitting" of (X).

\begin{lemma}  For each $X\in\{A,G,B\}$ we have 

i) $\mbox{Im}(r_X)=W_X$
(which justifies the uniqueness in the definition above.)

ii) $s_X$ is the composition 
$$
H^1_{dR}(X)\lra H^1_{dR}(X)/W_X\stackrel{i_X^{-1}}{\lra}\indif(X).
$$
\end{lemma}

{\bf Proof: } i) Let us recall a few simple facts about slope decomposition
of Frobenius modules (see for example \cite{zink} or \cite{berthelot}).
First, the unit root subspace is functorial with respect to Frobenius
morphisms, i.e. with respect to $K$-linear maps which commute with the 
Frobenii. Then the unit root subspace functor is left exact.
Let us now prove the lemma. Since $r_X$ is a section, it is an injective 
$K$-linear map. By assumption, $A$ has ordinary reduction, which implies that $B$ 
has
ordinary reduction, i.e. dim$_K(W_A)=$dim$_K(H(A))$ and 
dim$_K(W_B)=$dim$_K(H(B))$.
Therefore Im$(r_A)=W_A$ and Im$(r_B)=W_B$. Moreover, from the exact sequence
of Frobenius modules
$$
0\lra H^1_{dR}(B)\lra H^1_{dR}(G)\lra H^1_{dR}(T)\lra 0,
$$
we get an exact sequence of $K$-vector spaces
$$
0\lra W_B\lra W_G\lra (H^1_{dR}(T))^{slope=0}.
$$
As the last vector space is $0$, we get that dim$_K(W_G)=$dim$_K(W_B)$
which implies that Im$(r_G)=W_G.$

ii) Follows by definition from i).\hfill$\Box$
~

\begin{theorem}
Let $r_B: H^1(B,\cO) \rightarrow \hdrb$ be the unit root splitting on $B$. Then 
$L(r_B)$ is equal to $r_A$, the unit root splitting of $\hdra$. Hence $r_A$ is the 
unique splitting of the Hodge filtration of $\hdra$ making the following diagram 
commutative:
\begin{eqnarray*}
\begin{CD}
 @. \hdra  @<{r_A}<< \ha\\
 @. @VV{\pi^\ast}V @VV{\gamma}V \\
 (3)\quad \quad  @. \hdrg @. H\\
 @. @AA{p^\ast}A @VV{\beta\inv}V \\
@.\hdrb @<{r_B}<< {\hb}\\
\end{CD}
\end{eqnarray*}
 
\end{theorem}

{\bf Proof: } As $p^\ast$ is a morphism of Frobenius modules we have 
that $p^\ast (W_B)\subset W_G$. 
Therefore  $p^\ast r_B\beta^{-1}(H)
\subset W_G$ and as $p^\ast r_B\beta^{-1}$ is a section it follows that 
$p^\ast r_B\beta^{-1}=r_G$.

Now we want to prove that $\alpha^{-1} s_G\pi^\ast=s_A.$ For this let us 
remark that $\pi^\ast$ is a morphism of Frobenius modules, therefore
$\pi^\ast(W_A)\subset W_G$ hence  we have the following diagram:
$$
\begin{array}{cccccc}
\indif(A)&\stackrel{(i_A)^{-1}}{\lla}&H^1_{dR}(A)/W_A&
\lla &H^1_{dR}(A)\\
\alpha\downarrow&&\downarrow\pi^*&&\downarrow\pi^*\\
\indif(G)&\stackrel{(i_G)^{-1}}{\lla}&H^1_{dR}(G)/W_G&
\lla &H^1_{dR}(G)
\end{array}
$$
Both small squares are commutative therefore the large rectangle is 
commutative as well, so we have $\alpha s_A=s_G\pi^\ast$.

It remains to check uniqueness. Two splittings of the Hodge filtration 
which both make our diagram commutative differ by a homomorphism  
$\ha \rightarrow \mbox{\rm Inv}(A)$, which becomes zero after composition with 
$\pi^\ast$. Since $\pi^\ast$ is an isomorphism on invariant differentials our 
original map must also be zero.\hfill$\Box$

~

\section{The Mazur-Tate height corresponds to the unit root splitting}

Let $F$ be a rigid analytic group functor over $K$, i.e. a contravariant functor 
from the category  of rigid analytic $K$-varieties to the category of groups. We 
denote by $K[\epsilon]= K[T]/(T^2)$ the ring of dual numbers over $K$. Then the 
Lie algebra associated to $F$ is defined as
\[L(F)  = \mbox{\rm ker}(F(\mbox{Sp}\,\eps) \longrightarrow F(\mbox{Sp}\, K)) .\]

If $Z$ is a rigid analytic $K$ group variety, its Lie algebra is defined as the 
Lie algebra of the corresponding group functor, i.e. $\lie(Z) = \mbox{\rm ker} 
(Z(\eps) \rightarrow Z( K))$. We have a natural duality between $\indif(Z)$ and 
$\lie(Z)$.
Note that if $Z = X\an$ for some algebraic $K$-variety $X$ then we have a GAGA-
isomorphism $\lie(X\an) \simeq \lie \, X$, compatible with the dual isomorphism
$\indif(X\an) \simeq \indif(X)$.

Now we can define a rigidification of an exact sequence of rigid analytic 
$K$-groups as a splitting of the corresponding sequence of Lie algebras. Note that 
for an abelian variety $A$ over $K$ the GAGA-isomorphism
$\ext^1(A, \mathg_{m}) \simeq \ext^1(A\an, \mathg_{m}\an)$ can be continued to an 
isomorphism
$\mbox{\rm Extrig}(A, \mathg_{m}) \simeq \mbox{\rm Extrig}(A\an, \mathg_{m}\an)$. 
(When we use Ext groups for rigid analytic groups, we will always work in the site 
of all rigid $K$-varieties endowed with their Grothendieck topologies.)

Note further that the $K$-rational points of any smooth rigid analytic $K$-group  $Z$ 
form a Lie group over $K$ in the sense of \cite{bo}, III, \S 1. Its Lie algebra 
coincides with the Lie algebra of $Z$ as defined above. Rigid analytic morphisms 
of $K$-varieties give rise to analytic maps in the sense of \cite{bo2}, i.e. on 
$K$-rational points they are locally given by converging power series. (See also
\cite{se1}.) From now on 
we mean analytic in this sense when we talk about analytic maps. Whenever we  
consider rigid analytic objects, we will say so explicitely.

The uniformization maps $\pi$ and $\pi'$ which we introduced in the preceeding 
section induce isomorphisms on rigid analytic open subgroups
\[ G\an \supset \overline{G} \simeq \overline{A} \subset A\an\]
respectively
\[ G\dan \supset \overline{G}' \simeq \overline{A}' \subset A\dan. \]

By \cite{we}, 3.1, we have a (uniquely determined) isomorphism of biextensions
 of $G\an$ and 
$G\dan$ by $\mathg_{m}\an$
\[\theta: (\pi \times \pi')^\ast \pa\an {\longrightarrow} (p \times 
p')^\ast \pbb\an.\]

We fix a continuous, ramified homomorphism $\rho: K\tim \rightarrow \mathq_p$. 
Recall that by 
\cite{za}, p. 319, there is a non-zero $\Q_p$-linear map $\delta: K \rightarrow 
\Q_p$ and a branch $\lambda$ of the $p$-adic logarithm such that $\rho = \delta 
\circ \lambda$.

Let $J'$ be the universal vectorial extension of $B$, and let $\eta_B: B'(K) 
\rightarrow J'(K)$ be the unique splitting such that $\lie \,\eta_B = r_B$, the 
unit root splitting on $B$. By \cite{co}, Theorem 3.3.1, 
 the corresponding height pairing 
$\pair_{r_B}$ coincides with the Mazur-Tate pairing on $B$.

Now we consider the Mazur-Tate pairing on $A$. 
 Let $\tau_A$ denote the $\lambda$-splitting defined by the Mazur-Tate condition 
$\tau_A(x) = \lambda(\tilde{\sigma} (x))$ on $P^f_{\cA^0 \times \cA'}(R)$, where 
$\tilde{\sigma}:P^f_{\cA^0 \times \cA'} \rightarrow \mathg_m^\wedge $ is defined 
via the unique formal trivialization  of $P_{\cA^0 \times \cA'}^f$. Define 
$\tau_B$ in an analogous way. Then $\sigma_A = \delta \circ \tau_A$ and $\sigma_B 
= \delta \circ \tau_B$ are the $\rho$-splittings giving rise to the Mazur-Tate 
height pairings on $A$ respectively $B$. 

Fix some $a' \in \overline{A}'(K)$ and consider $x \in P_{\overline{A} \times 
\{a'\}}(K)$. The same argument as in \cite{we}, 7.1 shows that we have 
\[ \tau_B (pr \circ \theta(j(x))) = \tau_A(x)\]
where $j: P_{\overline{A} \times \{a'\} }\an \hookrightarrow  (\pi \times 
\pi')^\ast \pa\an$ is induced by $\overline{A}\simeq \overline{G} \hookrightarrow 
G$ and $\overline{A}'\simeq \overline{G}' \hookrightarrow G'$ and where  $pr: (p 
\times p')^\ast \pbb\an \rightarrow \pbb\an$ is the natural projection. 

This implies that the Mazur-Tate pairing on $A$ corresponds to the following 
section $\eta_A: A'(K) \rightarrow I'(K)$ of the projection: For $a'$ in the open 
subgroup $\ov{A}'$ denote its preimage in $\overline{G}'$ by $g'$, and put 
$b'=p'(g')$. Then $\theta$ induces an isomorphism $\pi^\ast P_{A \times \{a'\}}\an 
\simeq p^\ast P_{B \times \{b'\}}\an$. The map $\eta_B$ endows $ P_{B \times 
\{b'\}}$ with a rigidification. This induces a rigidification on $\pi^\ast P_{A 
\times \{a'\}}$, and hence on $P_{A \times \{a'\}}$, since $\pi$ is an isomorphism 
in a neighbourhood of the unit element. Then $\eta_A(a')$ is the extension $ P_{A 
\times \{a'\}}$ together with this rigidification. Since $I'$ is an extension of 
$A'$ by a vector group, 
$\eta_A$ can be extended uniquely to the whole of $A'(K)$.
Our next goal is to prove that $\eta_A$ is analytic as a map of $K$-varieties. 
As a first step we show the following lemma.

\begin{lemma} The $\lambda$-splitting $\tau_A: P(K) \rightarrow 
K$ which is induced by the formal splitting $P_{\cA^0 \times \cA'}^f \rightarrow 
\mathg^\wedge_m$ is an analytic map. 
\end{lemma}

{\bf Proof: }
Let us first show that under the  natural inclusion 
$P_{\cA^0\times \cA'}^f(R)$ may be identified with an open subset of $P(K)$ 
(in the Lie group topology as in \cite{bo2}.)
This follows from a more general fact, namely let 
$\cX$ be separated and of finite type over $R$ and let $\barX$ denote 
its special fiber. Let $Z\subset \barX$ be a closed subscheme. 
Denote by $\hX_Z$ the formal completion of $\cX$
along $Z$ and by $X=\cX_K$ the generic fiber of $\cX$. Let us denote 
by $X^{an}$ the rigid analytic variety corresponding to 
$X$ and by $\hX=\hX_{\barX}$ the formal completion of $\cX$ along
$\barX$.
Then, $\hX_Z$ is a formal scheme whose affine opens are formal spectra of 
$R$-algebras which are quotients of algebras of the form: $R\langle 
T_1,...T_r\rangle [[X_1,...,X_s]]$. To such a formal scheme one attaches
canonically (see \cite{berth}, 0.2) a  rigid analytic space, $(\hX_Z)_K$ together with a 
specialization map $sp: (\hX_Z)_K \rightarrow \hX_Z$ such that the following 
diagram is commutative
$$
\begin{array}{ccc}
(\hX_Z)_K&\stackrel{sp}{\longrightarrow}&\hX_Z\\
i\ \cap&&\downarrow\\
(\hX)_K&\stackrel{sp}{\longrightarrow}&\hX\\
j\ \cap\\
X^{an}
\end{array}
$$
by \cite{berth}, 0.2.7.
Here $(\hX)_K$ denotes the rigid analytic generic fiber of the formal scheme 
$\hX$ which is topologically of
finite type.
Moreover, $i((\hX_Z)_K)=sp^{-1}(Z)$ is an admissible open of $(\hX)_K$ 
and $j$ is an open immersion. Therefore, $(\hX_Z)_K$ can be identified with 
an admissible 
open of $X^{an}$, so that $(\hX_Z)_K(K)$ is an open subset of $X(K)$.
Let us now remark that $\hX_Z(R)=(\hX_Z)_K(K)$, which 
proves that $\hX_Z(R)$ is an open subset of $X(K)$.

In order to finish the proof of the lemma note that 
$P^f_{\cA^0\times \cA'}(R)$ may be identified with the set of points of $P(R)$
which project to points on $\cA^0(R)\times \cA'(R)$ whose specialization is 
$(0,0)$ on $\cA_k^0(k)\times \cA_k'(k)$. As $P^f_{\cA^0\times \cA'}$ is the 
trivial biextension, we fix formal 
sections $\tilde{\sigma}: P^f_{\cA^0\times \cA'}\lra {\mathbb{G}}_m^\wedge$ 
and 
$s: \cA^f\times \cA'^f \lra P^f_{\cA^0\times \cA'}$.

Let us now recall how the map $\tau_A :P(K)\lra K$ is defined.
Let $x\in P(K)$ and let $(a,a')\in A(K)\times A'(K)=\cA(R)\times \cA'(R)$
be its projection. For suitable integers $m$ and $n$ we have
 $(ma,na')\in \cA^f(R)\times \cA'^f(R)$, so let $y=s(ma,na')\in 
P^f_{\cA^0\times \cA'}(R)\subset P(K)$. Actually, $y$ is an analytic function 
of $x$. Then the elements $(m,n)x$ (biextension multiplication) and $y$ of $P(K)$ 
differ by a unique 
element $c\in K^*$, and let us note that $c$ is also analytic 
as a function of $x$ as the biextension operations are algebraic, hence 
analytic. Now we have
$$
\tau_A(x)=\frac{1}{mn} \left(\lambda(\tilde{\sigma}y)+\lambda(c)\right).
$$
Therefore
$\tau_A$ is analytic.\hfill$\Box$
~

We will now prove in general that analytic $\lambda$-splittings lead to analytic
splittings of the projection map $I'(K) \rightarrow A'(K)$. 

\begin{proposition}
Let $\tau: P(K) 
\rightarrow K$ be an analytic $\lambda$-splitting \`a la Mazur and Tate, where
$\lambda: K\tim \rightarrow K$ is a branch of the $p$-adic logarithm. 
Recall from section 1 that $\tau$ induces a splitting $\eta = \eta(\tau):
 A'(K) \rightarrow I'(K)$ by associating to $a'$ the extension $P_{A \times \{a'\}}$ 
 endowed with the rigidification  given by $\lie (\tau_{|P_{A \times \{a'\}}})$.
 Then $\eta$ is an analytic map.
\end{proposition}

{\bf Proof: } Choose Zariski open neighbourhoods $U \subset A$ and $U' \subset A'$ 
of the unit sections, such that the $\mathg_m$-torsor $P$ is trivial over $U 
\times U'$, i.e. we have a $\mathg_m$-equivariant isomorphism $\varphi: P_{U 
\times U'} \iso \mathg_{m} \times U \times U'$ over $U \times U'$. We can correct 
$\varphi$ by a suitable morphism $U' \rightarrow \mathg_{m}$ to achieve that 
$\varphi \circ e_{P/A'}: U' \rightarrow \mathg_{m } \times U \times U'$ is the map 
$u' \mapsto (1,1,u')$. Here $e_{P/A'}$ is the unit section of $P$ regarded as an 
$A'$-group. Then the map 
\[h:  U \times U' \stackrel{1 \times id}{\longrightarrow} \mathg_{m} \times U 
\times U' \stackrel{\varphi\inv}{\longrightarrow} P_{U \times U'}\]
is a $U'$-morphism compatible with the unit sections over $U'$. Therefore it 
induces a section $\Inf^1(U \times U' / U') \rightarrow \Inf^1(P_{U \times 
U'}/U')$ of the map  $\Inf^1(P_{U \times U'}/U') \rightarrow \Inf^1(U \times U' / 
U')$ given by the projection $P_{U \times U'} \rightarrow U \times U'$.  Here 
$\Inf^1$ denotes the first infinitesimal neighbourhood with respect to the unit
section. Now $U 
\times U'$ is an open subset of $A \times U'$ containing the image of the unit 
section over $U'$. Hence $\Inf^1(U \times U' /U') \simeq \Inf^1(A \times U'/U')$. 
Similarly, we have $\Inf^1(P_{U \times U'}/U') \simeq \Inf^1(P_{A \times U'}/U')$. 
Hence we get a commutative diagram
\begin{eqnarray*}
\begin{CD}
\Inf^1(A \times U'/U') @>>> \Inf^1(P_{A \times U'}/U')\\
 @VVV @VVV \\
 A \times U' @<<< P_{A \times U'},\\
\end{CD}
\end{eqnarray*}
which induces a rigidification on the extension $P_{A \times U'}$ of $A \times U'$ 
by $\mathg_{m,U'}$ over $U'$. 

In this way we get an element in $I'(U')$ projecting to the natural inclusion $U' 
\hookrightarrow A'$ in $A'(U')$, i.e. a local section $s: U' \rightarrow I'$ of 
the projection map $I'\rightarrow A'$. Since $s$ is a morphism of schemes, the 
corresponding map on $K$-rational points $s: U'(K) \rightarrow I'(K)$ is analytic. 

For every $x \in U'(K)$ its image $s(x)$ corresponds to the extension $P_{A \times 
\{x\}}$ together with the rigidification induced by the map $h(-,x): U 
\rightarrow P_{U \times \{x\}}$. 
It is easy to see that this rigidification can also be described as follows:
The isomorphism of $\mathg_m$-torsors $\varphi: P_{U \times U'} \rightarrow 
\mathg_{m} \times U \times U'$ induces a $K\tim$-equivariant analytic map $\psi: 
P_{U \times U'}(K) \rightarrow K\tim$, which in turn induces for all $x \in U'$ an 
analytic map $\psi_x: P_{U \times \{x\}}(K) \rightarrow K\tim$ respecting the unit 
sections. Hence the corresponding map on Lie algebras splits the Lie algebra 
sequence associated to $P_{A \times \{x\}}$, and gives the desired rigidification. 

Since $\tau: P(K) \rightarrow K$ is a $\lambda$-splitting, we
have $\tau(\alpha x) = 
\lambda(\alpha) + \tau(x)$ for all $\alpha \in K\tim$ and $x \in P(K)$. An 
analogous formula holds for $\lambda \circ \psi$. By hypothesis, the difference
$ \tau -\lambda \circ \psi: P_{U \times U'}(K) \rightarrow K$ is a 
$K\tim$-invariant analytic map, hence it factors over some analytic map $\theta: 
U(K)
 \times U'(K)
 \rightarrow K$. Since $\tau$ and $\lambda \circ\psi$ map the  unit section over 
$U'$ to zero, we have
$\theta(1,x) = 0$ for all $x \in U'(K)$. 

For all $x \in U(K)$ we now have two rigidifications on the extension 
$P_{A \times \{x\}}$. One comes from the point $s(x) \in I'(K)$. As shown above, 
it is given by
the map $\psi_x$. The other one comes from the point $\eta(x) \in I'(K)$, 
by definition it is given by the Lie algebra map corresponding to the 
Lie group homomorphism $\tau: P_{A \times \{x\}}(K) \rightarrow K$,
where we always identify $\lie(K\tim) \simeq \lie(K)$ by means of $\lie \lambda$.

A straightforward calculation shows that for all $x \in U'(K)$ these 
two  
rigidifications on $P_{A \times \{x\}}$ differ by the invariant 
differential $\omega_x$ corresponding to the Lie 
algebra map
\[ \lie \,\theta(-,x): \lie\, A \longrightarrow \lie \,K \simeq K.\]

This means that the sections $\eta: U'(K) \rightarrow I'(K)$ and 
$s: U'(K) \rightarrow I'(K)$ of the projection map differ by 
\begin{eqnarray*}
\omega: U'(K) & \longrightarrow & \indif(A)\\
 x  & \longmapsto  & \omega_x.
 \end{eqnarray*}
We claim that the map $\omega$ is analytic.
As the claim has a local nature let us fix $x_0\in U'(K)$ and choose  neighbourhoods
$V \subset U(K)$  of $1$ and $V' \subset U'(K)$ of $x_0$ such that 
$\theta$ is given by a convergent power series
\begin{eqnarray*}
\lefteqn{\theta (u_1,...,u_r;x_1,...,x_s)}\\
= & \sum_{i_1,...,i_r,j_1,...,j_s\ge 0}
a_{i_1,...,i_s;j_1,...,j_s}(u_1-1_1)^{i_1}...(u_r-1_r)^{i_r}(x_1-x_{0,1})^
{j_1}...(x_s-x_{0,s})^{j_s}
\end{eqnarray*}
on $V \times V'$. Here $(1_1,...,1_r)$ and $(x_{0,1},
...,x_{0,s})$ are the local coordinates of $1$ and $x_0$, respectively. A 
simple calculation shows that under the  identification 
$\mbox{Hom}_K(\lie(A), K))\simeq K^r$ given by the local coordinates
$u_1,\ldots, u_r$
we have
\begin{eqnarray*}
\lefteqn{
\lie(\theta)(-, x_1,...,x_s)} \\
= (&\sum_{j_1,...,j_s\ge 0}a_{1,0,...,0;j_1,...,j_s}
(x_1-x_{0,1})^{j_1}...(x_s-x_{0,s})^{j_s},...,\\
~ &\sum_{j_1,...,j_s\ge 0}
a_{0,...,0,1;j_1,...,j_s}(x_1-x_{0,1})^{j_1}...(x_s-x_{0,s})^{j_s}),
\end{eqnarray*}
which proves the claim.

Since the difference between $\eta$ and $s$ is an analytic map, we find that 
$\eta$ is analytic on $U'(K)$, hence everywhere.\hfill$\Box$
~

\begin{corollary}

The section $\eta_A: A'(K) \rightarrow I'(K)$ is analytic.

\end{corollary}

Therefore we know now that the Mazur-Tate $\lambda$-splitting $\tau_A$ is 
induced by a splitting of the Hodge filtration, namely by $\lie \, \eta_A$. 
Recall from the end of section 1 that this means that the Mazur-Tate height
pairing coincides with the $p$-adic height pairing induced by $\lie \, \eta_A$.
Hence it remains to identify $\lie \,\eta_A$ with the unit root splitting.
We need some technical preparations.

Let $X$ be a scheme over the base scheme $S$ or a rigid analytic variety over 
a rigid analytic base variety $S$. Then we denote by $\pic^\natural(X)=
\pic^\natural(X/S)$ the group of invertible sheaves on $X$ endowed with an 
integrable connection. For a morphism $X \rightarrow Y$ of $S$-schemes or rigid 
analytic varieties over $S$ we have a natural map $\pic^\natural(Y) \rightarrow 
\pic^\natural(X)$.

Let $(U_i)_i$ be a (Zariski or admissible rigid analytic) covering of $X$ such 
that the invertible sheaf $\cL$ is trivial on each $U_i$ with transition functions 
$f_{ij} \in \Gamma(U_i \cap U_j, \cO^\times)$. Then a connection $\nabla$ on $\cL$ 
gives rise to a collection of forms $\omega_i \in\Gamma(U_i, \Omega_{X/S}^1)$ such 
that $\omega_i - \omega_j = d f_{ij}/f_{ij}$ on $U_i \cap U_j$. Conversely, every 
such collection of forms defines a connection. $\nabla$ is integrable iff all the 
$\omega_i$ are closed forms. In this way one sees that
there is a functorial isomorphism
\[\pic^\natural(X) \simeq {\mathbb{H}}^1(X,\Omega^\times_{X/S}),\] 
where $\Omega^\times_{X/S}$ is the complex of sheaves
\[\Omega^\times_{X/S} = (\cO_X^\times \stackrel{d \log}{\longrightarrow} 
\Omega_{X/S}^1 \longrightarrow \Omega_{X/S}^2 \longrightarrow \ldots).\]

Now let $X$ be a smooth rigid analytic $K$-group. There is a natural 
functorial homomorphism
\[\mbox{\rm Extrig}_K(X, \mathg_{m}\an) \longrightarrow \pic^\natural(X)\]
defined as follows.
First of all, note that the rigidifications of the extension
\[0 \longrightarrow \mathg_{m}\an \longrightarrow Z \longrightarrow X 
\longrightarrow 0\]
correspond bijectively to invariant 1-forms $\omega$ on $Z$ which restrict to the 
form $dt/t$ on $\mathg_{m}\an$ (where $t$ is the standard parameter on 
$\mathg_{m}\an$), cf. \cite{za}, p.323.

Starting with an extension $Z$ and an invariant $1$-form $\eta$ on $Z$ restricting to 
$dt/t$ on $\bG_m^{an}$ we take a covering $(U_i)_i$ of $X$ such that the $\mathg_{m}\an$-torsor 
$Z$ is trivial over $U_i$, i.e. we have local sections $s_i: U_i \rightarrow Z$ of 
the projection $Z \rightarrow X$. 

Now we put $\omega_i= s_i^\ast \eta \in \Gamma(U_i, \Omega^1_{X/K})$. Let $\cL$ be 
the sheaf of sections of the line bundle associated to the torsor $Z$. Then, 
tautologically, the sections $f_{ij}\in \Gamma(U_i \cap U_j, \cO^\times)$ such 
that $f_{ij} s_j = s_i$ on $U_i \cap U_j$ are transition functions for $\cL$. One 
can easily check that the $\omega_i$ define a connection on $\cL$. All $\omega_i$ 
are closed since $\eta$ is closed as an invariant differential. Hence our 
connection is integrable.

\begin{lemma}
The following diagram commutes:
\begin{eqnarray*}
\begin{CD}
 \overline{A}'(K)  @>{\eta_A}>> \mbox{\rm Extrig}(A\an, \mathg_{m}\an) @>>> 
\pic^\natural(A\an)\\
 @AA{\simeq}A @VVV @VV{\pi^\ast}V\\
 \overline{G}'(K) @. \mbox{\rm Extrig}(G\an, \mathg_{m}\an) @>>> 
\pic^\natural(G\an)\\
@VVV @AAA @AA{p^\ast}A \\
B'(K) @>{\eta_B}>> \mbox{\rm Extrig}(B\an, \mathg_{m}\an) @>>> \pic^\natural 
(B\an).\\
\end{CD}
\end{eqnarray*}
\end{lemma}

{\bf Proof: }
The right part is commutative by functoriality, the left one by definition of 
$\eta_A$.\hfill$\Box$
~

\begin{lemma}
Let $X$ be a commutative rigid analytic group over $K$. Then we have a functorial 
isomorphism of $K$-vector spaces
\[L (\pic^\natural(X)) \iso H^1_{dR}(X),\]
where $L (\pic^\natural(X))$ is the Lie algebra of the group functor 
$S \mapsto \pic^\natural(X \times S)$.
\end{lemma}

{\bf Proof: }
This can be proved as in the algebraic case (see \cite{me}, 2.6.8).
Let $\Omega^\bullet$ be the complex $ (0 \rightarrow {\cal O}_X \rightarrow 
\Omega_{X/K} \rightarrow \Omega_{X/K}^2 \rightarrow \ldots)$ on $X$.  If we put 
$X_\epsilon = X \times_K \eps$ and denote by $p: X_\epsilon \rightarrow X$ the 
projection, there is a split exact 
sequence of complexes of abelian sheaves on $X$:
\[ 0 \longrightarrow \Omega^\bullet_{X/K} \longrightarrow  p_\ast 
\Omega_{X_\epsilon/\eps}^\times \longrightarrow \Omega_{X/K}^\times 
\longrightarrow 0.\]
Taking first hypercohomologies, our claim follows.\hfill$\Box$
~

Now we are able to prove the following theorem

\begin{theorem}
The Mazur-Tate height pairing on $A$ coincides with the height pairing defined
 by the unit root splitting $r_A$ of 
$\hdra$.
\end{theorem}
{\bf Proof: } It suffices to show that $r_A = \lie \, \eta_A$. 

Comparing $\eta_A$ with an algebraic splitting as in the proof of 3.2 we see that
we can pass from the diagram in Lemma 3.4 to the corresponding Lie algebra diagram.
Using Lemma 3.5, we  get a 
commutative diagram 
\begin{eqnarray*}
\begin{CD}
 \lie \,\overline{A}'(K)   @>{\mbox{\rm \small Lie} \, \eta_A}>> L 
(\pic^\natural(A\an)) @>{\simeq}>> \hdra\\
 @VV{\simeq}V @VVV @V{\pi^\ast}VV\\
 \lie \, \overline{G}'(K) @. L( \pic^\natural(G\an)) @>{\simeq}>> \hdrg\\
@VVV @AAA @A{p^\ast}AA \\
\lie \, B'(K) @>>> L (\pic^\natural(B\an)) @>{\simeq}>> \hdrb.\\
\end{CD}
\end{eqnarray*}
Here the upper horizontal map is a section of the natural projection $\hdra 
\rightarrow \ha \simeq \lie\, \overline{A}'(K)$ and the lower horizontal map is 
induced by $\lie \,\eta_B = r_B$. Besides, the diagram
\begin{eqnarray*}
\begin{CD}
\ha   @>{\simeq}>> \lie \,\overline{A}'(K)\\
 @VVV @VV{\simeq}V\\
H @. \lie \,\overline{G}'(K) \\
@VVV @VVV \\
\hb @>{\simeq}>> \lie \, B'(K)\\
\end{CD}
\end{eqnarray*}
commutes by \cite{ls2}, 6.7. 
Hence $\lie \,\eta_A$ makes diagram $(3)$ in Theorem 2.2 commutative, and our 
claim follows from this result.\hfill$\Box$
~

\section{Appendix}

In this section we will prove that the Hodge filtration of the first de Rham 
cohomology group of a semiabelian variety over $K$ is given by its invariant 
differentials. This answers a question raised in  \cite{ls2}.

Let us first recall some notations from section 2. Let $K$ be a finite 
extension of ${\Bbb Q}_p$, $T=({\Bbb G}_m)^t$ a split torus over $K$,
$B$ an abelian variety over $K$ with good reduction and $G$
a semiabelian variety over $K$ defined by the extension of algebraic 
groups over $K$:
$$
(4)\quad 0\lra T\stackrel{g}{\lra} G \stackrel{p}{\lra} B\lra 0.
$$
We will prove

\begin{theorem} For $X\in \{T,G,B\}$ the image of $\indif(X)$ in
$H^1_{dR}(X)$ can be naturally identified with the Hodge filtration on
$H^1_{dR}(X)$.
\end{theorem}

Note that this result is well known for $X=B$.

As $T$ and $G$ are smooth schemes over $K$ which are not proper, let us 
first briefly review the Hodge theory for non-proper varieties.

Let $X$ be a smooth scheme over $K$, $Y$ a smooth and proper scheme over $K$ and 
$X\stackrel{i_X}{\lra}Y$ an open embedding. Let $Z=Y-X$ and suppose that 
$Z$ is a divisor with normal crossings over $K$. Denote by $\Omega^i_Y(\log Z)$
the sheaf of $i$-th differential forms on $Y$ with logarithmic poles 
along $Z$. The first deRham cohomology of $Y$ with logarithmic poles along $Z$ is 
defined as $H^i_{dR}(Y,\log(Z)) = {\Bbb H}^i(\Omega_Y^\bullet(\log Z))$.
Then we have (see \cite{deligne}):

1) The inclusion $\Omega_Y^\bullet (\log Z)\subset (i_X)_*\Omega_X^\bullet$
induces an isomorphism
$$
H^i_{dR}(X)\cong H^i_{dR}(Y,\log(Z)).$$

2) The spectral sequence
$$
E_1^{i,j}=H^j(Y, \Omega^i_Y(\log Z)) ==> H^{i+j}_{dR}(Y, \log Z)=
H^{i+j}_{dR}(X)
$$
degenerates at $E_1$ and induces the Hodge filtration.

3) The Hodge filtration thus defined on $H^i_{dR}(X)$ does not depend on the 
compactification of $X$.

As a consequence of the above we have that the first step of the Hodge filtration 
on $H^1_{dR}(X)$, denoted $F^0_X$, is given by
$$
F^0_X=\mbox{Ker}(d:H^0(Y, \Omega_Y^1(\log Z))\lra H^0(Y, \Omega^2_Y(\log Z))).
$$

Let us now prove the theorem.

{\bf Proof: } a) We will first prove the statement for $X=T$. Consider the
 embedding $T=({\Bbb G}_m)^t\stackrel{i_T}{\lra} \PP$, where each $\bG_m$ is 
embedded naturally in ${\Bbb P}^1_K$, after a parameter was chosen. 
This embedding may be seen as a smooth compactification of $T$. Moreover
the group law: $T\times T\longrightarrow T$ extends naturally to a 
morphism $T\times \PP\longrightarrow \PP$, i.e. endows $\PP$ with a natural
action of $T$. (See \cite{serre}.)

 Let $Z:=\PP-T$. It is clearly a divisor with normal crossings and 
so, by the above we have $H^1_{dR}(T)=H^1_{dR}(\PP, \log Z)$. As the 
$1$-forms in $\indif(T)$ are closed and are logarithmic when 
considered in $\bP:=\PP$, we have an inclusion 
$$
\indif(T)\lra F^0_T=\mbox{Ker}(d:H^0(\bP, \Omega^1_\bP(\log Z)\lra 
H^0(\bP,\Omega^2_\bP(\log Z))).
$$
To show that the map is surjective let us fix parameters $z_1,z_2,...,z_t$
on each factor of $\bP=\PP$ compatible with the parameters on the factors 
of $T=(\bG_m)^t$. Then an element  $\omega\in H^0(\bP, \Omega^1_\bP(\log Z))$
has the form 
$$
\omega=\sum_{n=1}^t f_n\frac{dz_n}{z_n},
$$
where $f_n\in H^0(\bP, \cO_\bP)=K$. Therefore $\omega$ is closed and 
$\omega\in \indif(T)$.

b) Let us now prove the theorem for $X=G$.
We will first briefly recall how one obtains a good compactification of $G$
(see \cite{serre}). Let's recall the short exact sequence (4) at the 
beginning of this section. As torus torsors are locally trivial in the 
Zarski topology, we may find a finite open, affine cover $\{U_i\}_i$ of $B$
such that 
$$
p^{-1}(U_i)\cong U_i\times T\stackrel{1 \times i_T}{\lra}U_i\times \PP.
$$
Now embed each $U_i\times T\subset U_i\times \PP$
and notice that the isomorphisms used to glue the opens $\{U_i\times T\}_i$
extend naturally and define gluing data for the set $\{U_i\times \PP\}_i$.
We glue the schemes $\{U_i\times \PP\}_i$ along these isomorphisms 
and obtain a proper scheme $M$ over $K$ such that the following hold

i) The natural map $G\stackrel{i_G}{\lra} M$ is an open embedding.

ii) We have a natural morphism $M\lra B$ and the torsor
 structure of $G$ over $B$ naturally extends to a stucture of
principal fiber-space of $M$ over $B$ of fiber type $\PP$.

iii) Let $V_i:=U_i\times \PP\subset M$. Then if $Z:=M-G$ we have $Z|_{V_i}=
U_i\times (\PP-T)=U\times Z_i$ which is a divisor with normal crossings in 
$V_i$ for every $i$. Therefore $Z$ is a divisor with normal crossings in $M$.

As a consequence, we have $H^1_{dR}(G)=H^1_{dR}(M, \log Z)$ and as the 
invariant differentials on $G$ are closed we have a natural inclusion
$$
\indif(G)\lra \mbox{Ker}(d: (i_G)_*\Omega^1_{G}(M)\lra 
(i_G)_*\Omega^2_{G}(M)).
$$
Let us fix $\{Y_j\}_j$ a standard open affine cover of $\PP$ (compatible
with the parameters $z_1,z_2,...,z_t$ chosen above.) Then 
$\{U_i\times Y_j\}_{i,j}$ is an open affine cover of $M$ and let 
$Z_{ij}:=Z\cap (U_i\times Y_j)=U_i\times (Y_j-T)$. Let $\omega\in \indif
(G)$. As both $\Omega_M^1(\log Z)$ and $(i_G)_*\Omega^1_{G}$ are coherent 
sheaves on $M$, in order to show that $\omega\in \Omega_M^1(\log Z)(M)$
it would be enough to show that 
$$
\omega|_{U_i\times Y_j}\in \Omega^1_M(\log Z)(U_i\times Y_j)=\Omega^1_{U_i\times 
Y_j}(\log Z_{ij}).
$$
Let us fix $i,j$ and denote by $p_1,p_2$ the projections from $U_i\times Y_j$
to its factors (in this order). Then we have a natural isomorphism
$$
\Omega^1_{U_i\times Y_j}(\log Z_{ij})\cong p_1^*\Omega^1_{U_i}\times p_2^*
\Omega_{Y_j}(\log (Y_j-T)).
$$
Moreover,
$$
\omega|_{U_i\times Y_j}\in ((i_G)_*\Omega^1_{G})(U_i\times Y_j)=
\Omega^1_{U_i\times (Y_j\cap T)}=p_1^*\Omega^1_{U_i}\times 
p_2^*\Omega^1_{Y_j\cap T}.
$$
So $\omega|_{U_i\times Y_j}=(\omega_1, \omega_2)$, with $\omega_1\in 
p_1^*\Omega_{U_i}^1$ and $\omega_2=p_2^*((g^*\omega)|_{Y_j\cap T}).$
As $\omega\in \indif(G)$ it follows that $g^*\omega\in \indif(T)$ and we 
have proved at a) above that $g^*\omega$ is logarithmic, i.e. that 
$g^*\omega|_{Y_j\cap T}\in \Omega^1_{Y_j}(\log (Y_j-T))$.

This shows that we have a natural injection:
$$
\indif(G)\lra F^0_G.
$$
Now we claim that the exact sequence $(4)$ induces an exact sequence
$$
(5)\quad 0\lra F^0_B\lra F^0_G\lra F^0_T\lra 0.
$$
To see this let us recall that for all $X\in \{T,G,B\}$ we have isomorphisms
as Gal$(\overline{K}/K)$-modules 
$$
H^1_{et}(X_{\overline{K}}, {\Bbb Q}_p(1))\cong V_p(X),
$$
where $V_p(X):=(\stackrel{\mbox{lim}}{\leftarrow} 
X_K[p^n](\overline{K}))\otimes_{{\Bbb Z}_p}
{\Bbb Q}_p$ is the Tate module of $X_K$. The exact sequence $(4)$ also 
produces an exact sequence of Gal$(\overline{K}/K)$-modules (see \cite{ls2})
$$
0\lra V_p(B)\lra V_p(G)\lra V_p(T)\lra 0.
$$
Using Fontaine's theory and the results in \cite{tsuji}
we have a commutative diagram of {\it filtered modules} with exact rows:
$$
\begin{array}{ccccccccc}
0&\lra&D_{dR}(V_p(B))&\lra&D_{dR}(V_p(G))&\lra&D_{dR}(V_p(T))&\lra&0\\
&&\downarrow\cong&&\downarrow\cong&&\downarrow\cong\\
0&\lra&H^1_{dR}(B)&\lra&H^1_{dR}(G)&\lra&H^1_{dR}(T)&\lra&0
\end{array}
$$
This implies the exactness of the sequence $(5)$ which implies that
dim$_K(F^0_G)=$dim$_K(\indif(G))$ (we have used that the theorem is true 
for $T$ and $B$). Therefore the natural injection
$\indif(G)\lra F^0_G$ is an isomorphism.
~

~\\[4ex]
\parbox{40ex} {Adrian Iovita \\ University of Washington \\
Dept. of Mathematics  \\ Box 354350 \\ Seattle, WA 98195 - 4350 \\ email: 
iovita@math.washington.edu}
\parbox{40ex}{\begin{flushright} Annette Werner\\Universit\"at 
M\"unster\\Mathematisches Institut\\
Einsteinstr. 62 \\ D - 48149 M\"unster \\ email: werner@math.uni-muenster.de
\end{flushright}}

\end{document}